\documentclass[12pt,a4paper]{article}
\usepackage[left=2cm,right=2cm,  top=2cm,bottom=2cm,bindingoffset=0cm]{geometry}

\usepackage{color}

\usepackage{bm}
\usepackage{dcolumn}
\usepackage{subfig}
\usepackage{flushend}
\usepackage{indentfirst}
\usepackage{graphics}
\usepackage{amsmath}
\usepackage{graphicx}
\usepackage{epstopdf}
\usepackage{float}
\usepackage{amssymb}
\usepackage[font=footnotesize,labelfont=bf]{caption}
\usepackage[labelsep=period]{caption}

\newcommand{\be}{\begin{equation}}
\newcommand{\ee}{\end{equation}}

\usepackage{latexsym,amsmath}
\usepackage{hyperref}
\usepackage{amsfonts}

\newtheorem{thm}{Theorem}[section]

\newtheorem{cor}{Corollary}

\newtheorem{defn}{Definition}[section]

\begin{document}

\title{Estimation of periodically correlated random fields that are isotropic on a sphere}

\date{}

\maketitle

\author{\textbf{Iryna Golichenko}$^{1}$, \textbf{Oleksandr Masyutka}$^{2}$, \textbf{Mykhailo Moklyachuk}$^{3*}$, \\\\
$^{1}$\emph{ Department of Mathematical Analysis
and Probability Theory, National Technical University of Ukraine
``Kyiv Polytechnic Institute'' Kyiv 03056, Ukraine}\\
$^{2}$\emph{ Department of Mathematics and Theoretical Radiophysics,
Taras Shevchenko National University of Kyiv, Kyiv 01601, Ukraine}\\
$^{3}$\emph{ Department of Probability Theory, Statistics and Actuarial
Mathematics, Taras Shevchenko National University of Kyiv, Kyiv 01601, Ukraine}\\
$^{*}$Corresponding Author: Moklyachuk@gmail.com}\\\\\\

\noindent \textbf{Abstract.} \hspace{2pt}

The problem of optimal linear estimation of functionals
depending on the unknown values of a spatial temporal isotropic random field
$\zeta(j,x)$, which is periodically correlated
 with respect to discrete time argument $j\in\mathrm Z$ and mean-square continuous isotropic on the unit sphere ${S_n}$ with respect to spatial argument $x\in{S_n}$.
Estimates are based on
observations of the field $\zeta(j,x)+\theta(j,x)$ at points $(j,x):$ $j\in Z\backslash\{0, 1, .... , N\}$,  $x\in S_{n}$,
 where $\theta(j,x)$ is an
uncorrelated with $\zeta(t,x)$
spatial temporal isotropic random field, which is periodically correlated
 with respect to discrete time argument $j\in\mathrm Z$ and mean-square continuous isotropic on the sphere ${S_n}$ with respect to spatial argument $x\in{S_n}$.
Formulas for calculating the mean square errors and the spectral characteristics of the optimal linear
estimate of the functional are derived in the case where the spectral density matrices are exactly known. Formulas that determine the least
favourable spectral density matrices and the minimax (robust) spectral
characteristics are proposed in the case where the spectral density matrices are not exactly known but a class of admissible spectral density matrices is given.\\

\noindent \textbf{Keywords:} \hspace{2pt} spatial temporal isotropic random field, isotropic random field, periodically correlated random field, robust estimate, mean
square error, least favourable spectral density, minimax spectral
characteristic.\\

\noindent \textbf{AMS subject classifications:} 60G60; 62M40; 62M20; 93E10; 93E11

\section{Introduction}

Cosmological Principle (first formulated by Einstein): the Universe is,
in the large, homogeneous and isotropic \cite{Bartlett}. Last
decades indicate growing interest to the
spatio-temporal data measured on the surface of a sphere. These data
includes cosmic microwave background (CMB) anisotropies \cite{Bartlett}, \cite{Hu}, \cite{Kogo}, \cite{Okamoto}, \cite{Adshead}, medical imaging \cite{Kakarala}, global and land-based temperature data \cite{Jones}, \cite{SubbaRao1}, gravitational and geomagnetic data, climate model \cite{North}. Some basic results and references on the theory of isotropic random fields on a sphere can be found in the books \cite{Yadrenko} and \cite{Yaglom} (1987). For more recent applications and results see books\cite{Gaetan}, \cite{Cressie}, \cite{Marinucci} and several
papers covering a number of problems in general for spatial temporal isotropic
observations \cite{SubbaRao2}, \cite{Terdik}.

Periodically correlated processes and fields are not homogeneous but have numerous properties similar to properties of
stationary processes and fields. They describe appropriate models of numerous physical
and man-made processes. A comprehensive list of the existing
references up to the year 2005 on periodically correlated processes
and their applications was proposed in \cite{Serpedin}. See also reviews \cite{Antoni}. For more details see a survey paper \cite{Gardner} and book \cite{Hurd}. Note, that in the literature periodically correlated processes are named in multiple different ways such as cyclostationary, periodically nonstationary or cyclic correlated processes.

Among the current trends of the theory of stochastic processes and fields important is the direction which focuses on the problem of estimation of unknown values of random processes and fields.
The problem of estimation of random processes and fields includes interpolation, extrapolation and filtering problems.

The mean square optimal estimation problems for periodically
correlated with respect to time isotropic on a sphere random fields
are natural generalization of the linear extrapolation,
interpolation and filtering problems for stationary stochastic
processes and homogeneous random fields.
Effective methods of solution of the linear extrapolation, interpolation and filtering problems for stationary stochastic processes were developed under the condition of certainty where spectral densities of processes
and fields are known exactly (see, for example, \cite{Kolmogorov}, survey article \cite{Kailath} by Kailath, books \cite{Rozanov}, \cite{Wiener}, \cite{Yaglom}, \cite{Yadrenko}, articles \cite{Moklyachuk:1979}, \cite{Moklyachuk:1980}.

Particularly relevant in recent years is the problem of estimation of values of processes and fields under uncertainty where the spectral densities of processes
and fields are not known exactly.
Such problems arise when considering problems of automatic control theory, coding and signal processing in radar and sonar, pattern recognition problems of speech signals and images. See a survey paper \cite{Kassam}.

The classical approach to the problems of interpolation, extrapolation and filtering of stochastic processes and random fields is based on the assumption that the spectral densities of processes and fields are known. In practice, however, complete information about the spectral density is impossible in most cases. To overcome this complication one finds parametric or nonparametric estimates of the unknown spectral densities or select these densities by other reasoning. Then applies the classical estimation method provided that the estimated or selected density is the true one. This procedure can result in a significant increasing of the value of error as Vastola and Poor \cite{Vastola} have demonstrated with the help of some examples. This is a reason to search estimates which are optimal for all densities from a certain class of admissible spectral densities. These estimates are called minimax since they minimize the maximal value of the error of estimates.
A survey of results in minimax (robust) methods of data processing can be found in the paper \cite{Kassam}.
 The paper by Grenander \cite{Grenander} should be marked as the first one where the minimax approach to extrapolation problem for stationary processes was
developed. For more details see, for example, books \cite{Moklyachuk:2008}, \cite{Moklyachuk:2012}, \cite{Golichenko}.
In papers \cite{Dubovetska12}, \cite{Dubovetska13}, \cite{Dubovetska14_1}, \cite{Dubovetska14_2} the minimax-robust estimation problems (extrapolation, interpolation and filtering) are investigated for the linear functionals which depend on unknown values
of periodically correlated stochastic processes. Methods of solution the minimax-robust estimation problems for time-homogeneous isotropic random fields on a sphere were developed in \cite{Moklyachuk:1994}, \cite{Moklyachuk:1995}, \cite{Moklyachuk:1996}.
In papers \cite{Dubovetska14_M}, \cite{Dubovetska15_1}, \cite{Dubovetska15_2} results of investigation of minimax-robust estimation problems for periodically correlated isotropic random fields are described.

In this article we deal with the problem of optimal linear
estimation of the functional
\[A_{N}\zeta=\sum _{j=0}^{N}\int_{S_{n}}a(j,x)\zeta(j,x)m_{n}(dx)\]
depending on the unknown values of a spatial temporal isotropic random field
$\zeta(j,x)$, which is periodically correlated
(cyclostationary with period $T$) with respect to discrete time argument $j\in\mathrm Z$ and mean-square continuous isotropic on the unit sphere ${S_n}$ in Euclidean space ${\mathbb E}^n$ with respect to spatial argument $x\in{S_n}$.
Estimates are based on
observations of the spatial temporal isotropic random field $\zeta(j,x)+\theta(j,x)$ at points
$(j,x):$ $j\in Z\backslash\{0, 1,\dots,N\}$; $x\in{S_n}$, where $\theta(j,x)$ is an
uncorrelated with $\zeta(t,x)$
spatial temporal isotropic random field, which is periodically correlated
(cyclostationary with period $T$) with respect to discrete time argument $j\in\mathrm Z$ and mean-square continuous isotropic on the unit sphere ${S_n}$ with respect to spatial argument $x\in{S_n}$.
Formulas are
derived for computing the value of the mean-square error and the
spectral characteristic of the optimal linear estimate of the
functional $A_N\zeta$ in the case of spectral certainty, where the spectral density matrices of the
fields are known. Formulas are proposed that determine the least
favourable spectral density matrices and the minimax-robust spectral
characteristic of the optimal estimate of the functional $A_N\zeta$
in the case of spectral uncertainty, where the
spectral density matrices are not known exactly, but, instead, some special classes $D =D_f \times D_g$ of
admissible spectral density matrices are given.
\\
We use the Kolmogorov (see, for example, \cite{Kolmogorov}) Hilbert space projection method based on properties of the Fourier coefficients of some functions from the spectral density matrices.

\section*{Spectral properties of periodically correlated random fields}
Let $S_{n}$ be a unit sphere in the $n$-dimensional Euclidean space ${\mathbb E}^n$, let $m_{n}(dx)$ be the Lebesgue measure on $S_{n}$, and let
\[S_{m}^{l}(x),\,\, l=1, \dots, h(m,n);\,\, m=0, 1,\dots\]
be the orthonormal spherical harmonics of degree $m$ \cite{Muller}.

A spatial temporal isotropic random field ${\zeta(j,x)}$, $j\in\mathbb Z$, $x\in{S_n}$ is called periodically correlated
(cyclostationary with period $T$) with respect to discrete time argument $j\in\mathrm Z$ and mean-square continuous isotropic on the sphere ${S_n}$ with respect to spatial argument $x\in{S_n}$ if
\[
{\mathbb E}{\zeta(j,x)}=0,\quad {\mathbb
E}|{\zeta(j,x)}|^2<\infty,
\]
\[ {\mathbb E}\left({\zeta(j+T,x)} \overline
{\zeta (k+T,y)}\right)=B\left(j,k,\cos\langle x,y\rangle \right),
\]
where $\cos\langle x,y\rangle=(x,y)$ is a ``distance'' between points $x,y\in{S_n}$.
The spatial temporal isotropic random field ${\zeta(j,x)}$ can be represented in the form
\[ {\zeta(j,x)}=
{\sum_{m=0}^{\infty}} {\sum_{l=1}^{h(m,n)}} S_m^l (x)\zeta_{m}^l
(j),
\]
\[ \zeta_{m}^l (j)= {\int_{S_n}}{\zeta(j,x)} S_m^l
(x)\,m_n(dx),
\]
where
\[\zeta_{m}^l (j),\,\, j\in\mathbb Z, \,\,m=0,1,\dots;\,\, l=1,\ldots,h(m,n)\]
are mutually uncorrelated periodically correlated stochastic sequences with the correlation functions $b^{\zeta}_m(j,k)$:
\[
{\mathbb E}\left(\zeta_{m}^l (j+T)\overline{\zeta_{u}^v
(k+T)}\right)=\delta_m^u \delta_l^v\,\,b^{\zeta}_m(j,k),\]
\[ m,u=0,1,\dots;\,\,
l,v=1,\ldots,h(m,n);\,j,k\in\mathbb Z.
\]
The correlation function of the field $\zeta(t,x)$ can be represented as follows
\[
B\left(j,k,\cos\langle x,y\rangle \right)=\]
\[=\frac{1}{\omega_n}\sum_{m=0}^{\infty}h(m,n) \frac
{C_{m}^{(n-2)/2}(\cos\langle x,y\rangle)} {C_{m}^{(n-2)/2}(1)}\,\,
b^{\zeta}_m(j,k),
\]
where $\omega_n=(2\pi)^{n/2}\Gamma(n/2)$, and  $C_{m}^{(n-2)/2}(z)$, $m=0,1,\dots,$ are the Gegenbauer polynomials \cite{Muller}.

It follows from \cite{Gladyshev} (see also \cite{Makagon}) that the stochastic sequence $\zeta_{m}^l (j)$, $j\in\mathbb Z$, is periodically correlated with period $T$ if and only if there exists a $T$-variate stationary sequence
\[ \vec\xi_m^l(j)= \{\xi_{mk}^l(j)\}_{k=0}^{T-1},\quad
j\in\mathbb Z,\]
such that $\zeta_{m}^l (j)$, $j\in\mathbb Z$, can be represented in the form

\[
\zeta_{m}^l (j)= \sum_{k=0}^{T-1} e^{2\pi ijk/T}\xi_{mk}^l(j),\quad
j\in\mathbb Z.
\]
The sequence $\vec\xi_m^l(j)=\{\xi_{mk}^l(j)\}_{k=0}^{T-1}$,  $j\in\mathbb Z$, is called generating sequence of the periodically correlated sequence $\zeta_{m}^l (j)$,  $j\in\mathbb Z$.

Denote by $\Phi_m^{\vec{\xi}}(d\lambda)$ the matrix spectral measure function of the $T$-variable vector stationary sequence $
\vec\xi_m^l(j)= \{\xi_{mk}^l(j)\}_{k=0}^{T-1}$ resulting from the Gladyshev representation. Denote by $\Phi_m^{\vec\zeta}(d\lambda)$ the matrix spectral measure function of the $T$-variable vector stationary sequence
\[\vec\zeta_m^l(j)= \{\zeta_{mk}^l(j)\}_{k=0}^{T-1},\quad [\vec{\zeta}_m^l(j)]_{k}=\zeta_m^l(jT+k)
\]
\[
j\in {\mathbb
Z},\quad k=0,1,\dots,T-1,
\]
arising from the splitting into blocks of length $T$ the univariate periodically correlated sequence $\zeta_m^l(t)$. The relation of spectral matrices $\Phi_m^{\vec{\xi}}(d\lambda)$ and $\Phi_m^{\vec{\zeta}}(d\lambda)$ is described by the formula
\[
\Phi_m^{\vec{\zeta}}(d\lambda)=T\cdot V(\lambda)\Phi_m^{\vec
{\xi}}(d\lambda/T)V^{-1}(\lambda),
\]
where $V(\lambda)$ is an unitary $T\times T$ matrix whose $(k, j)$-th element is of the form
\[
v_{kj}(\lambda)=\frac{1}{\sqrt{T}}e^{2{\pi} ikj/T+ik\lambda/T},\quad k,j=
0,1,\dots,T-1.
\]
This relation can also be expressed as
\[
\Phi_m^{\vec
{\xi}}(d\lambda)=\frac{1}{T}\cdot V^{-1}(T\lambda) \Phi_m^{\vec{\zeta}}
(Td\lambda)V(T\lambda).
\]
Consequently, if there exists the spectral density matrix $F_{m}^{\vec {\xi}}(\lambda)$ of the $T$-variate stationary sequence $\vec {\xi}_m^l(j)$ then there exists the spectral density matrix $F_{m}^{\vec {\zeta}}(\lambda)$ of the $T$-variate stationary sequence $\vec {\zeta}_m^l(j)$ and these two density matrices satisfy the relation
\[
F_{m}^{\vec {\zeta}}(\lambda)=T\cdot V(\lambda)F_{m}^{\vec
{\xi}}(\lambda/T)V^{-1}(\lambda).
\]

\section*{Hilbert space projection method of estimation}

Consider the problem of the mean square optimal linear estimation of the functional
\[
A_{N}\zeta ={\sum_{j=0}^{N}}\,\,\,{\int_{S_n}} \,\,a(j,x)\zeta
(j,x)\,m_n(dx)
\]
depending on the unknown values of a spatial temporal isotropic random field
$\zeta(j,x)$, which is periodically correlated
 with respect to discrete time argument $j\in\mathrm Z$ and mean-square continuous isotropic on the unit sphere ${S_n}$ in Euclidean space ${\mathbb E}^n$ with respect to spatial argument $x\in{S_n}$.
 Estimates are based on observations of the spatial temporal isotropic random field $\zeta(j,x)+\theta(j,x)$ at points $(j,x):$ $j\in Z\backslash\{0, 1, .... , N\}$, $x\in S_{n}$, where $\theta(j,x)$ is an uncorrelated with $\zeta(j,x)$ spatial temporal isotropic random field, which is periodically correlated
 with respect to discrete time argument $j\in\mathrm Z$ and mean-square continuous isotropic on the sphere ${S_n}$ with respect to spatial argument $x\in{S_n}$. The field $\theta(j,x)$ has the representation
 \[ {\theta(j,x)}=
{\sum_{m=0}^{\infty}} {\sum_{l=1}^{h(m,n)}}S_m^l (x)  \theta_{m}^l
(j)=\]
\[={\sum_{m=0}^{\infty}} {\sum_{l=1}^{h(m,n)}} S_m^l (x)
\sum_{k=0}^{T-1} e^{2\pi ikj/T} \eta_{mk}^l (j),
\]
\[ \theta_{m}^l (j)= {\int_{S_n}}{\theta(j,x)} S_m^l (x)\,m_n(dx).
\]
In this representation
\[\theta_{m}^l (j), j\in\mathbb
Z, m=0,1,\dots;\,\, l=1,\ldots,h(m,n)\]
are mutually uncorrelated periodically correlated stochastic sequences with the correlation functions $b^{\theta}_m(j,s)$:
\[
{\mathbb E}\left(\theta_{m}^l (j+T)\overline{\theta_{u}^v
(k+T)}\right)=\delta_m^u \delta_l^v\,\,b^{\theta}_m(j,k),\]
\[ m,u=0,1,\dots;\,\,
l,v=1,\ldots,h(m,n);\,j,k\in\mathbb Z,
\]
and $\vec\eta_m^l(j)= \{\eta_{mk}^l(j)\}_{k=0}^{T-1}$ are vector-valued stationary sequences generating the periodically correlated sequences $\theta_{m}^l (j)$.

We will suppose that there exist the matrix spectral densities $F_{m}(\lambda)=F_{m}^{\xi}(\lambda)$ and $G_{m}(\lambda)=G_{m}^{\eta}(\lambda)$ of the generating sequences
$\vec\xi_m^l(j)= \{\xi_{mk}^l(j)\}_{k=0}^{T-1}$ and $\vec\eta_m^l(j)= \{\eta_{mk}^l(j)\}_{k=0}^{T-1}$
and denote by $M(F+G)$ the set of those $m\in\mathbb Z$ for which
\begin{equation}\label{eq1}
\int _{-\pi }^{\pi }\, Tr\, \left[\left(F_{m}(\lambda )+G_{m}(\lambda )\right)^{-1} \right]\, d\lambda \, <\, \infty .
\end{equation}
We will consider random fields which satisfy the following minimality condition: $M(F+G)\not=\emptyset$.
Under this minimality condition the errorless estimation of the functional is impossible \cite{Rozanov}.

Making use the Gladyshev results \cite{Gladyshev} we can represent the functional $A_{N}\zeta$ in the form
\[
A_{N}\zeta = {\sum_{m=0}^{\infty}}
{\sum_{l=1}^{h(m,n)}}{\sum_{j=0}^{N}} (\vec{a}_m^l(j))^{\top}\vec{\xi}_m^l (j),
\]
\[
\vec{a}_m^l(j)=({a}_{m0}^l(j),{a}_{m1}^l(j),\dots,{a}_{m(T-1)}^l(j))^{\top},
\]
\[
{a}_{mk}^l(j)={a}_m^l(j) e^{2\pi ikj/T},\,k=0,1,\dots,T-1,
\]
where $\vec\xi_m^l(j)= \{\xi_{mk}^l(j)\}_{k=0}^{T-1}$  are vector-valued stationary sequences generating the periodically correlated sequences $\zeta_{m}^l (j)$.

Every linear estimate $\widehat{A_{N}\zeta}$ of the functional ${A_{N}\zeta}$ is determined by the spectral stochastic measures
\[(Z_{\vec{\xi}})_m^l(d\lambda)=\left\{(Z_{\vec{\xi}})_{mk}^l(d\lambda)\right\}_{k=0}^{T-1},\quad
(Z_{\vec{\eta}})_m^l(d\lambda)=\left\{(Z_{\vec{\eta}})_{mk}^l(d\lambda)\right\}_{k=0}^{T-1},
\]
of the generating sequences $\vec\xi_m^l(j)=
\{\xi_{mk}^l(j)\}_{k=0}^{T-1}$ and $\vec\eta_m^l(j)=
\{\eta_{mk}^l(j)\}_{k=0}^{T-1}$,
 and by the spectral characteristic
\[h(\lambda)=\left\{h_{m}^{l}(\lambda): l=1, ... , h(m,n), m=0, 1, ...\right\} \]
of the estimate $\widehat{A_{N}\zeta}$ which belongs to the space $L_{2}^{N-} (F+G)$ generating by the functions
\[ h_{m}^{l}(\lambda)=\sum _{j\in Z\backslash\{0, 1, .... , N\}}\vec{h}_{m}^{l}(j)e^{ij\lambda },\]
\[l=1, ... , h(m,n),\; m=0, 1, ...,\]
which satisfy the condition
\[\sum _{m=0}^{\infty }\sum _{l=1}^{h(m,n)}\int _{-\pi }^{\pi }h_{m}^{l}(\lambda )^{\top}H_{m}(\lambda )\overline{h_{m}^{l}(\lambda )}d\lambda <\infty,\]
where $H_{m}(\lambda)=F_{m}(\lambda)+G_{m}(\lambda)$.

\noindent The estimate is of the form
\[\widehat{A\zeta}
= \sum_{m\in M} {\sum_{l=1}^{h(m,n)}}\int_{-\pi}^{\pi}
h_{m}^{l}(\lambda )^{\top}\left(\lambda\right) \left(
(Z_{\vec{\xi}})_m^l(d\lambda) +
(Z_{\vec{\eta}})_m^l(d\lambda)\right).
\]

\noindent The mean square error $\Delta (h;F,G)={\mathbf E}\left|A_{N}\zeta-\widehat{A_{N}\zeta}\right|^{2}$ of the estimate $\widehat{A_N\zeta}$
is determined by matrices of spectral densities
\[F(\lambda)=\{F_{m}(\lambda):m=0,1\dots\},\quad
G(\lambda)=\{G_{m}(\lambda):m=0,1\dots\}\]

\noindent of the generating sequences $\vec\xi_m^l(j)=
\{\xi_{mk}^l(j)\}_{k=0}^{T-1}$ and $\vec\eta_m^l(j)=
\{\eta_{mk}^l(j)\}_{k=0}^{T-1}$, and the spectral characteristic
$h\left(\lambda\right)$ of the estimate.

\noindent The value of the mean square error can be calculated by the formula
\[\Delta (h;F,G)={\mathbf E}\left|A_{N}\zeta-\widehat{A_{N}\zeta}\right|^{2}=\]
\[=\sum_{m\in M}\sum _{l=1}^{h(m,n)}\left\{\frac{1}{2\pi } \int _{-\pi }^{\pi }(A_{mN}^{l}(\lambda)-h_{m}^{l}(\lambda))^{\top} F_{m}(\lambda )\overline{(A_{mN}^{l}(\lambda)-h_{m}^{l}(\lambda)}d\lambda\right\}+
\]
\[
+\sum_{m\in M}\sum _{l=1}^{h(m,n)}\left\{\frac{1}{2\pi } \int _{-\pi }^{\pi }h_{m}^{l}(\lambda)^{\top} G_{m}(\lambda )\overline{h_{m}^{l}(\lambda)}d\lambda\right\},\]
\[A_{mN}^{l}(\lambda) =\sum _{j=0}^{N}\vec{a}_{m}^{l}(j)e^{ij\lambda }.\]
The spectral characteristic $h(F,G)$ of the optimal linear estimate $\widehat{A_{N}\zeta}$ minimizes the mean square error.
\[
\Delta(F,G)=\Delta(h(F,G);F,G)=\]
\[=\mathop {\min }\limits_{\widehat
{A_{N}\zeta}}{\mathbf E}\left|A_{N}\zeta-\widehat{A_{N}\zeta}\right|^{2}.
\]

\noindent In this article we use the Kolmogorov \cite{Kolmogorov} Hilbert space projection method based on properties of the Fourier coefficients of some functions from the spectral density matrices.
With the help of the proposed method we can derive formulas for calculation the mean square error $\Delta(F,G)=\Delta(h(F,G);F,G)$ and the spectral characteristic $h(F,G)$ of the optimal linear estimate $\widehat{A_{N}\zeta}$ of the functional ${A_{N}\zeta}$ under the condition that
spectral density matrices $F_{m}(\lambda )$ and $G_{m}(\lambda )$ are exactly known and satisfy the minimality condition: $M(F+G)\not=\emptyset$.
Following the method we found the optimal linear estimate $\widehat{A_{N}\zeta}$ as projection of ${A_{N}\zeta}$ on the closed linear subspace $H^{N-}_{\zeta+\theta}$ generated in the space $H=L_2(\Omega,\mathcal F,P)$ by values
\[\{\zeta(j,x)+\theta(j,x): j\in Z\backslash\{0, 1, .... , N\}; x\in S_n\}.\]
This projection is determined by conditions:
\[1)\,\widehat{A_{N}\zeta}\in H^{N-}_{\zeta+\theta};\;
2)\, (A_{N}{\zeta}-\widehat{A_{N}\zeta})\perp H^{N-}_{\zeta+\theta}.\]
The second condition is satisfied if for all $m=0,$ $1,\dots$; $l=1,\ldots,h(m,n)$, $j\in Z\backslash\{0, 1, \ldots , N\}:$
\[
\frac{1}{2\pi}{\int_{-\pi}^\pi}\left[(A_{mN}^{l}(\lambda)-h_{m}^{l}(\lambda))^{\top}F_{m}(\lambda)-h_{m}^{l}(\lambda)^{\top}G_{m}(\lambda)\right]e^{-ij\lambda}d\lambda=0.
\]
These relations mean that for all $m=0,1,\dots$; $l=1,\ldots,h(m,n)$:
\[
A_{mN}^{l}(\lambda)^{\top}F_{m}(\lambda)-h_{m}^{l}(\lambda)^{\top}H_{m}(\lambda)=
{C}_{mN}^{l}(\lambda)^{\top},
\]
\[
C_{mN}^{l}(\lambda) =\sum _{j=0}^{N}\vec{c}_{m}^{l}(j)e^{ij\lambda }.
\]
where $\vec{c}_{m}^{l}(j)$, $j=0,1,\dots, N$ are unknown coefficients. It follows from the indicated relations that the
spectral characteristic $h(F,G)$ of the optimal linear estimate$\widehat{A_{N}\zeta}$ is of the form
\[
h_{m}^{l}(F,G)^{\top} =r_{F}^{N}(\lambda)(H_{m}(\lambda ))^{-1}=
\]
\begin{equation}\label{eq2}
=A_{mN}^{l}(\lambda)^{\top}-r_{G}^{N}(\lambda)(H_{m}(\lambda ))^{-1},
\end{equation}
where
\[r_{F}^{N}(\lambda)=A_{mN}^{l}(\lambda)^{\top} F_{m}(\lambda )-C_{mN}^{l}(\lambda) ^{\top},\]
\[r_{G}^{N}(\lambda)= A_{mN}^{l}(\lambda)^{\top} G_{m}(\lambda )+C_{mN}^{l}(\lambda)^{\top}.\]
The first condition $\widehat{A_{N}\zeta}\in H^{N-}_{\zeta+\theta}$, is satisfied if for all $m=0,1,\dots$;
$l=1,\ldots,h(m,n)$ and $k=0,1,\dots, N$:
\[
\sum_{j=0}^{N}\left\{\frac{1}{2\pi}{\int_{-\pi}^\pi}[F_{m}(\lambda)(H_{m}(\lambda))^{-1}]^{\top}
e^{i(j-k)\lambda}d\lambda\right\}\vec{a}_m^l(j)=
\]
\begin{equation}\label{eq3}
=\sum_{j=0}^{N}\left\{\frac{1}{2\pi}{\int_{-\pi}^\pi}
[ (H_{m}(\lambda))^{-1}]^{\top}e^{i(j-k)\lambda}d\lambda\right\}\vec{c}_m^l(j).
\end{equation}
Let us introduce operators $\mathbf{B}_{mN}$, $\mathbf{D}_{mN}$,
$\mathbf{R}_{mN}$  determined by matrices $B_{mN}=\{B_{m}(k,j)\}_{k,j=0}^{N}$, $D_{mN}=\{D_{m}(k,j)\}_{k,j=0}^{N}$, $R_{mN}=\{R_{m}(k,j)\}_{k,j=0}^{N}$ composed with the help of the Fourier coefficients
\[B_{m} (k, j)=\frac{1}{2\pi } \int _{-\pi }^{\pi }\left[(H_{m}(\lambda ))^{-1} \right]^{\top} e^{i(j-k)\lambda } d\lambda ,\]
\[D_{m} (k, j)=\frac{1}{2\pi } \int _{-\pi }^{\pi }\left[F_{m}(\lambda ) (H_{m}(\lambda ))^{-1} \right]^{\top} e^{i(j-k)\lambda } d\lambda ,\]
\[R_{m} (k, j)=\frac{1}{2\pi } \int _{-\pi }^{\pi }  \left[F_{m}(\lambda ) (H_{m}(\lambda ))^{-1} G_{m}(\lambda )\right]^{\top} e^{i(j-k)\lambda } d\lambda, \] \[ k , j=0, 1, \ldots , N,\]
and vectors
\[\mathbf{a}_{mN}^{l} =\left\{\vec{a}_{m}^{l}(0),\ldots, \vec{a}_{m}^{l}(N)\right\},\]
\[\mathbf{c}_{mN}^{l} =\left\{\vec{c}_{m}^{l}(0), \ldots , \vec{c}_{m}^{l}(N)\right\}.\]
Taking into consideration the introduced operators and vectors we can represent equations (\ref{eq3}) in the form
\[ \mathbf{D}_{mN}\mathbf{a}_{mN}^{l}=\mathbf{B}_{mN}\mathbf{c}_{mN}^{l},\]
\[ m=0,1,\dots;\; l=1,\ldots,h(m,n).\]
This means that the unknown coefficients $\mathbf{c}_{mN}^{l}$ are determined by the equation
\[\mathbf{c}_{mN}^{l}= (\mathbf{B}_{mN})^{-1} \mathbf{D}_{mN}\mathbf{a}_{mN}^{l},\]
\[m=0,1,\dots;\;l=1,\ldots,h(m,n).\]
It follows from the derived relations that the value of the mean square error $ \Delta (F, G)$  of the
optimal linear estimate $\widehat{A_{N}\zeta}$ of the functional ${A_{N}\zeta}$ can be calculated by formula
\[
\Delta (h;F,G)=\]
\[=\sum _{m\in M}\sum _{l=1}^{h(m,n)}\frac{1}{2\pi } \int _{-\pi }^{\pi }\left[r_{G}^{N}(\lambda)(H_{m}(\lambda ))^{-1}F_{m}(\lambda )(H_{m}(\lambda ))^{-1}(r_{G}^{N}(\lambda))^{*}\right]d\lambda
+\]
\[
+\sum _{m\in M}\sum _{l=1}^{h(m,n)}\frac{1}{2\pi } \int _{-\pi }^{\pi }\left[r_{F}^{N}(\lambda)(H_{m}(\lambda ))^{-1}
 G_{m}(\lambda )(H_{m}(\lambda ))^{-1}(r_{F}^{N}(\lambda))^{*}\right]d\lambda=\]
\begin{equation}\label{eq4}
=\sum _{m\in M}\sum _{l=1}^{h(m,n)}\biggl(\left\langle \mathbf{c}_{mN}^{l}, \mathbf{B}_{mN} \mathbf{c}_{mN}^{l} \right\rangle
+\left\langle \mathbf{a}_{mN}^{l},  \mathbf{R}_{mN} \mathbf{a}_{mN}^{l} \right\rangle\biggr).
\end{equation}
Let us summarize our results and present them in the form of statements.

\begin{thm}
 Let $\zeta(j,x)$ and $\theta(j,x)$
be uncorrelated spatial temporal isotropic random fields, which are periodically correlated
 with respect to time argument $j\in\mathrm Z$ and isotropic on the sphere ${S_n}$ with respect to spatial argument $x\in{S_n}$
 which have the spectral density matrices $F_{m}(\lambda )$ and $G_{m}(\lambda )$ that satisfy
the minimality condition: $M(F+G)\not=\emptyset$.
 The spectral characteristic $h(F,G)$ and the value of the mean square error $\Delta (h;F,G)$  of the optimal linear estimate of the functional $A_N\zeta$ based on observations of the field $\zeta(j,x)+\theta(j,x)$ at points $(j,x):$ $j\in Z\backslash\{0, 1,\dots, N\}$, $x\in{S_n},$ are calculated by formulas (\ref{eq2}), (\ref{eq4}).
\end{thm}

\begin{cor}
Let $\zeta(j,x)$ be spatial temporal isotropic random field, which is periodically correlated
 with respect to time argument $j\in\mathrm Z$ and isotropic on the sphere ${S_n}$ with respect to spatial argument $x\in{S_n}$
 which have the spectral density matrix $F_{m}(\lambda )$ that satisfy
the minimality condition
\begin{equation}\label{eq5}
\int _{-\pi }^{\pi }\, Tr\, \left[\left(F_{m}(\lambda )\right)^{-1} \right]\, d\lambda  \, <\, \infty .
\end{equation}
The value of the mean square error $\Delta (F)$ and the spectral characteristic $h(F)$ of the optimal linear estimate of the functional $A_{N}\zeta$ based on observations of the field $\zeta(j,x)$ at points $(j,x):$ $j\in Z\backslash\{0, 1,\dots, N\}$, $x\in{S_n},$ can be calculated by formulas
\[
\Delta (F)=\sum _{m\in M}\sum _{l=1}^{h(m,n)}\frac{1}{2\pi } \int _{-\pi }^{\pi }\left[(C_{mN}^{l}(\lambda) )^{\top} )(F_{m}(\lambda ))^{-1}(C_{mN}^{l}(\lambda) )^{\top} )^{*}\right]d\lambda=\]
\begin{equation}\label{eq6}
=\sum _{m\in M}\sum _{l=1}^{h(m,n)}\left\langle \mathbf{B}_{mN}^{-1}\mathbf{a}_{mN}^{l}, \mathbf{a}_{mN}^{l} \right\rangle,
\end{equation}
\begin{equation} \label{eq7}
h_{m}^{l}(F)^{\top}=(A_{mN}^{l}(\lambda))^{\top}-(C_{mN}^{l}(\lambda))^{\top}(F_{m}(\lambda ))^{-1},
\end{equation}
where
\[ C_{mN}^{l}(\lambda) =\sum _{j=0}^{N}\vec{c}_{m}^{l}(j)e^{ij\lambda },\; \mathbf{c}_{mN}^{l}=\mathbf{B}_{mN}^{-1}  \mathbf{a}_{mN}^{l},\]
$\mathbf{B}_{mN}$ is operator determined by the matrix $B_{mN}=\{B_{m}(k,j)\}_{k,j=0}^{N}$ composed with the help of the Fourier coefficients
\[B_{m} (k, j)=\frac{1}{2\pi } \int _{-\pi }^{\pi }\left[(F_{m}(\lambda ))^{-1} \right]^{\top} e^{i(j-k)\lambda } d\lambda ,\]
\[k , j=0, 1,\dots, N.\]
\end{cor}

\section*{Minimax-robust method of interpolation} The proposed formulas may be employed under the condition that the spectral densities $F_{m}(\lambda )$,  $G_{m}(\lambda )$ of the fields $\zeta(j,x), \theta(j,x)$ are exactly known. In the case where the densities are not known exactly but a set $D=D_{F} \times D_{G}$ of possible spectral densities is given the minimax (robust) approach to estimation of functionals of the unknown values of random fields is reasonable. Instead of searching an estimate that is optimal for a given spectral densities we find an estimate that minimizes the mean square error for all spectral densities from given class simultaneously.
\begin{defn}
 The spectral densities $F^{0} (\lambda )$, $G^{0} (\lambda )$ are called least favourable in a given class $D=D_{F} \times D_{G}$ for the optimal linear interpolation of the functional $A_{N}\zeta$ if the following relation holds true
\[\Delta (h(F^{0} ,G^{0} );F^{0} ,G^{0} )=\mathop{\max }\limits_{(F,G)\in D} \Delta (h(F,G);F,G).\]
\end{defn}

It follows from the relationships (\ref{eq1}) -- (\ref{eq7}) that the following theorems hold true.
\begin{thm}
Spectral densities $F_{m}^{0} \left(\lambda \right), G_{m}^{0} \left(\lambda \right)$ that satisfy the minimality condition \eqref{eq1} are least favourable in the class  $D=D_{F}\times D_{G}$ for the optimal linear interpolation of the functional $A_{N}\zeta$ if they determine operators $\mathbf{B}_{mN}^{0} $, $\mathbf{D}_{mN}^{0} $, $\mathbf{R}_{mN}^{0}$ giving a solution to the extremum problem
\[
\mathop{\max }\limits_{\left(F, G\right)\in D} \sum _{m=0}^{\infty }\sum _{l=1}^{h(m,n)}\biggl(\left\langle \mathbf{B}_{mN}^{-1} \mathbf{D}_{mN} \mathbf{a}_{m}^{l}, \mathbf{D}_{mN} \mathbf{a}_{m}^{l} \right\rangle+\left\langle \mathbf{a}_{m}^{l}, \mathbf{R}_{mN} \mathbf{a}_{m}^{l} \right\rangle\biggr)=
\]
\begin{equation}\label{eq8}
=\sum _{m=0}^{\infty }\sum _{l=1}^{h(m,n)}\biggl(\left\langle (\mathbf{B}_{mN}^{0} )^{-1} \mathbf{D}_{mN}^{0} \mathbf{a}_{m}^{l}, \mathbf{D}_{mN}^{0} \mathbf{a}_{m}^{l} \right\rangle+\left\langle \mathbf{a}_{m}^{l} , \mathbf{R}_{mN}^{0} \mathbf{a}_{m}^{l} \right\rangle\biggr).
\end{equation}
\end{thm}

\begin{thm}
The spectral density $F_{m}^{0} \left(\lambda \right)$ that satisfies the minimality condition (\ref{eq5}) is least favourable in the class  $D_{F}$ for the optimal linear interpolation of the functional $A_{N}\zeta$  based on observations of the field $\zeta(j,x)$ at points $j\in Z\backslash\{0, 1,\dots, N\}$,  $x\in S_{n}$ if it determines operators $\mathbf{B}_{mN}^{0} $ giving a solution to the extremum problem
\begin{equation}\label{eq9}
\mathop{\max }\limits_{F\in D_{F}} \sum _{m=0}^{\infty }\sum _{l=1}^{h(m,n)}\left\langle \mathbf{B}_{mN}^{-1} \mathbf{a}_{m}^{l}, \mathbf{a}_{m}^{l} \right\rangle=\sum _{m=0}^{\infty }\sum _{l=1}^{h(m,n)}\left\langle (\mathbf{B}_{mN}^{0} )^{-1}\mathbf{a}_{m}^{l}, \mathbf{a}_{m}^{l} \right\rangle.
\end{equation}
\end{thm}

\begin{defn}
The spectral characteristic
\[h^{0} (\lambda )=\left\{h_{m}^{0l}(\lambda): l=1,\dots, h(m,n), m=0, 1, \dots\right\}\]
of the optimal linear interpolation of the functional $A_{N}\zeta$ is called minimax-robust if there are satisfied conditions
\[h^{0} (\lambda)\in H_{D} =\bigcap _{(F,G)\in D}L_{2}^{N-} (F+G) ,\]
\[\mathop{\min }\limits_{h\in H_{D} } \mathop{\sup }\limits_{(F,G)\in D} \Delta (h;F,G)=\mathop{\sup }\limits_{(F,G)\in D} \Delta (h^{0} ;F,G).\]
\end{defn}

The least favourable spectral densities  $F_{m}^{0} (\lambda )$, $G_{m}^{0} (\lambda )$ and the minimax (robust) spectral characteristic $h^{0} (\lambda)\in H_{D}$ form a saddle point of the function $\Delta (h;F,G)$. The saddle point inequalities
\[\Delta (h;F^{0} ,G^{0} )\geq \Delta (h^{0} ;F^{0} ,G^{0} )\ge \Delta (h^{0} ;F,G),\]
\[\forall (F,G)\in D,\; \forall h\in H_{D} \]
hold true if $h^{0} =h(F^{0} ,G^{0} )\in H_{D} $ and $(F_{m}^{0} ,G_{m}^{0} )$ is a solution to the conditional extremum problem
\begin{equation}\label{eq10}
\Delta (h(F^{0} ,G^{0} );F^{0} ,G^{0} )=\mathop{\sup }\limits_{(F,G)\in D} \Delta (h(F^{0} ,G^{0} );F,G),
\end{equation}
where
\[
\Delta (h(F^{0} ,G^{0} );F,G)=\]
\[=\sum _{m\in M}\sum _{l=1}^{h(m,n)}\frac{1}{2\pi } \int _{-\pi }^{\pi }\left[r_{G}^{N0}(\lambda)(H_{m}^{0}(\lambda ))^{-1}F_{m}(\lambda )(H_{m}^{0}(\lambda ))^{-1}(r_{G}^{N0}(\lambda))^{*}\right]d\lambda+
\]
\[
+\sum _{m\in M}\sum _{l=1}^{h(m,n)}\frac{1}{2\pi } \int _{-\pi }^{\pi }\left[r_{F}^{N0}(\lambda)(H_{m}^{0}(\lambda ))^{-1} G_{m}(\lambda )(H_{m}^{0}(\lambda ))^{-1}(r_{}^{N0}(\lambda))^{*}\right]d\lambda.
\]

The following theorems hold true.
\begin{thm}
Let $(F_{m}^{0} ,G_{m}^{0} )$ be a solution to the extremum problem (\ref{eq10}). The spectral densities $F_{m}^{0}$ and $G_{m}^{0}$ are  least favourable in the class $D=D_{F} \times D_{G}$, and the spectral characteristic $h^{0}=h(F^{0}, G^{0})$ is minimax spectral characteristic for the optimal interpolation of the functional $A_{N}\zeta$ if the condition $h(F^{0}, G^{0})\in H_{D}$ holds true.
\end{thm}

\begin{thm}
The spectral density $F^{0}(\lambda)\in D_{F}$ which satisfies the minimality condition (\ref{eq5}) is least favourable in the class $D_{F}$ for the optimal interpolation of the functional $A_{N}\zeta$ based on observations of the field $\zeta(j,x)$ at points $j\in Z\backslash\{0, 1, .... , N\}$,  $x\in S_{n}$ if $F^{0}(\lambda)$ is a solution to the extremum problem
\begin{equation}\label{eq11}
\Delta (h(F^{0});F^{0})=\mathop{\sup }\limits_{F\in D_{F}} \Delta (h(F^{0});F)
\end{equation}
where
\[
\Delta (h(F^{0});F)=\]
\[=\sum _{m\in M }\sum _{l=1}^{h(m,n)}\frac{1}{2\pi } \int _{-\pi }^{\pi }\left[(C_{mN}^{l0}(\lambda) )^{\top}(F_{m}^{0}(\lambda ))^{-1} F_{m}(\lambda )(F_{m}^{0}(\lambda ))^{-1}(C_{mN}^{l0}(\lambda) )^{\top} )^{*}\right]d\lambda.
\]
The spectral characteristic $h^{0}=h(F^{0})$ is the minimax spectral characteristic for the optimal interpolation of the functional $A_{N}\zeta$ if the condition $h(F^{0})\in H_{D}$ holds true.
\end{thm}

The conditional extremum problem(10) is equivalent to the following unconditional extremum problem
\begin{gather}\label{eq12}
\Delta _{D} (F,G)=-\Delta (h(F^{0} ,G^{0} );F,G)+\delta ((F,G)|D)\to \inf ,
\end{gather}
where $\delta ((F,G)|D)$ is the indicator function of the set D. A solution to the problem (\ref{eq12}) is characterized by the condition $0\in \partial \Delta _{D} (F^{0} ,G^{0} )$, where $\partial \Delta _{D} (F^{0} ,G^{0} )$ is the subdifferential of the convex functional $\Delta _{D} (F,G)$ at point $(F^{0} ,G^{0} )$.

\section*{Least favourable spectral densities in the class $D_{0}\times D_{\varepsilon}$}

Consider the problem of minimax estimation of the functional $A_{N}\zeta$
depending on the unknown values of a spatial temporal isotropic random field
$\zeta(j,x)$, which is periodically correlated
 with respect to discrete time argument $j\in\mathrm Z$ and isotropic on the sphere ${S_n}$ with respect to spatial argument $x\in{S_n}$
 based on observations of the spatial temporal isotropic random field
$\zeta(j,x)+\theta(j,x)$ at points $(j,x):$ $ j=0,-1,-2,\dots; x\in{S_n},$
under the condition that matrices of spectral densities
$F(\lambda)=\{F_m(\lambda):m=0,1\dots\}$ and
$G(\lambda)=\{G_m(\lambda):m=0,1\dots\}$ of the field $\zeta(j,x)$
and the field $\theta (j,x)$ are not known exactly, but the
following pairs of sets of spectral densities that give restrictions on the first moment and describe describe the model of "$\varepsilon$-contamination"  of spectral densities. The first pair is
\[D_{0}^{1} =\biggl\{F(\lambda )|\frac{1}{2\pi\omega_{n} }\sum_{m=0}^{\infty}h(m,n) \int _{-\pi }^{\pi }Tr F_{m}(\lambda )d\lambda =p \biggr\},\]
\[D_{\varepsilon}^{1} =\biggl\{G(\lambda )|Tr G_{m}(\lambda )=(1-\varepsilon)Tr U_{m} (\lambda )+\varepsilon Tr V_{m}(\lambda),\]
 \[
 \frac{1}{2\pi\omega_{n} }\sum_{m=0}^{\infty}h(m,n) \int _{-\pi }^{\pi } Tr G_{m}(\lambda )d\lambda  =q,\biggr\},\]

The second pair of sets of admissible spectral densities is
\[D_{0}^{2} =\biggl\{F(\lambda )|\frac{1}{2\pi \omega_{n}}\sum_{m=0}^{\infty}h(m,n) \int _{-\pi }^{\pi }f_{m}^{kk} (\lambda )d\lambda =p_{k} , k=\overline{1,T} \biggr\},\]
\[D_{\varepsilon}^{2} =\biggl\{G(\lambda )|g_{m}^{kk} (\lambda )=(1-\varepsilon)u_{m}^{kk} (\lambda )+\varepsilon v_{m}^{kk} (\lambda ), \]
\[\frac{1}{2\pi\omega_{n} }\sum_{m=0}^{\infty}h(m,n) \int _{-\pi }^{\pi }g_{m}^{kk} (\lambda )d\lambda  =q_{k} , k=\overline{1,T}\biggr\};\]

The third pair of sets of admissible spectral densities is
\[D_{0}^{3} =\biggl\{F(\lambda )|\frac{1}{2\pi \omega_{n}}\sum_{m=0}^{\infty}h(m,n) \int _{-\pi }^{\pi }\left\langle B_{1}, F_{m}(\lambda )\right\rangle d\lambda  =p\biggr\},\]
\[D_{\varepsilon}^{3} =\biggl\{G(\lambda )|\left\langle B_{2} , G_{m}(\lambda )\right\rangle=(1-\varepsilon)\left\langle B_{2} , U_{m} (\lambda )\right\rangle+\varepsilon \left\langle B_{2}, V_{m} (\lambda )\right\rangle,\]
\[
 \frac{1}{2\pi\omega_{n} }\sum_{m=0}^{\infty}h(m,n) \int _{-\pi }^{\pi }\left\langle B_{2} , G_{m}(\lambda )\right\rangle d\lambda  =q\biggr\};\]
The forth pair of sets of admissible spectral densities is
\[D_{0}^{4} =\biggl\{F(\lambda )|\, \frac{1}{2\pi \omega_{n}}\sum_{m=0}^{\infty}h(m,n) \int _{-\pi }^{\pi }F_{m}(\lambda )d\lambda =P \biggr\},\]
\[D_{\varepsilon}^{4} =\biggl\{G(\lambda )|G_{m}(\lambda )=(1-\varepsilon)U_{m} (\lambda )+\varepsilon V_{m} (\lambda ),\]
 \[ \frac{1}{2\pi \omega_{n}}\sum_{m=0}^{\infty}h(m,n) \int _{-\pi }^{\pi }G_{m}(\lambda )d\lambda  =Q\biggr\},\]

Here $V_{m} (\lambda ), U_{m} (\lambda )$ are given matrices of spectral densities, $p, q, p_{k}, q_{k}, k=\overline{1,T}$ are given numbers, $B_{1}, B_{2}, P, Q$ are given positive-definite Hermitian matrices. From the condition $0\in \partial \Delta _{D} (F^{0} ,G^{0} )$ we find equations which determine the least favourable spectral densities for these given pairs of sets of admissible spectral densities.

For the first pair $D_{0}^{1}\times D_{\varepsilon}^{1}$ we have equations
\begin{equation}\label{eq13}
\sum_{l=1}^{h(m,n)}(r_{G}^{N0}(\lambda))^{*}r_{G}^{N0}(\lambda)=\alpha_{m} ^{2} (H_{m} ^{0} (\lambda ))^{2},
\end{equation}

\begin{equation}\label{eq14}
\sum_{l=1}^{h(m,n)}(r_{F}^{N0}(\lambda))^{*}r_{F}^{N0}(\lambda)=
(\beta_{m}^{2} +\gamma _{m} (\lambda ))(H_{m} ^{0} (\lambda ))^{2}.
\end{equation}

For the second pair $D_{0}^{2}\times D_{\varepsilon}^{2}$ we have equations
\begin{equation}\label{eq15}
\sum_{l=1}^{h(m,n)}(r_{G}^{N0}(\lambda))^{*}r_{G}^{N0}(\lambda)= H_{m} ^{0} (\lambda )\left\{\alpha _{mk}^{2} \delta _{kn} \right\}_{k,n=1}^{T} H_{m} ^{0} (\lambda ),
\end{equation}

\begin{equation}\label{eq16}
\sum_{l=1}^{h(m,n)}(r_{F}^{N0}(\lambda))^{*}r_{F}^{N0}(\lambda)=H_{m} ^{0} (\lambda )
\left\{(\beta _{mk}^{2} +\gamma _{mk} (\lambda ))\delta _{kn} \right\}_{k,n=1}^{T} H_{m} ^{0} (\lambda ).
\end{equation}

For the third pair $D_{0}^{3}\times D_{\varepsilon}^{3}$ we have equations

\begin{equation}\label{eq17}
\sum_{l=1}^{h(m,n)}(r_{G}^{N0}(\lambda))^{*}r_{G}^{N0}(\lambda)=\alpha_{m} ^{2} H_{m} ^{0} (\lambda )B_{1}^{\top}  H_{m} ^{0} (\lambda ),
\end{equation}

\begin{equation}\label{eq18}
\sum_{l=1}^{h(m,n)}(r_{F}^{N0}(\lambda))^{*}r_{F}^{N0}(\lambda)=(\beta_{m}^{2} +\gamma_{m}^{'}(\lambda)) H_{m} ^{0} (\lambda )B_{2}^{\top} H_{m} ^{0} (\lambda ).
\end{equation}

For the fourth pair $D_{0}^{4}\times D_{\varepsilon}^{4}$ we have equations
\begin{equation}\label{eq19}
\sum_{l=1}^{h(m,n)}(r_{G}^{N0}(\lambda))^{*}r_{G}^{N0}(\lambda)= H_{m} ^{0} (\lambda ) \vec{\alpha_{m}}\cdot \vec{\alpha_{m}}^{*} H_{m} ^{0} (\lambda ),
\end{equation}
\begin{equation}\label{eq20}
\sum_{l=1}^{h(m,n)}(r_{F}^{N0}(\lambda))^{*}r_{F}^{N0}(\lambda)=H_{m} ^{0} (\lambda ) (\vec{\beta_{m}}\cdot \vec{\beta_{m}}^{*}+\Gamma _{m} (\lambda ))H_{m} ^{0} (\lambda ).
\end{equation}
Here $\Gamma _{m} (\lambda )$ are Hermitian matrices and
\[\gamma _{m} (\lambda )\le 0, \gamma _{mk} (\lambda )\le 0, k=\overline{1, T}\, a.e.,\]
\[\gamma_{m}^{'}(\lambda)\le 0, \Gamma _{m_{1}} (\lambda )\le 0,\, a.e.,\]
\[\gamma _{m} (\lambda )=0: Tr G_{m}^{0} (\lambda )>(1-\varepsilon)Tr U_{m}(\lambda ),\]
\[\gamma _{mk} (\lambda )=0: g^{0kk}_{m} (\lambda )>(1-\varepsilon)u^{kk}_{m}(\lambda ),\]
\[\gamma _{m}^{'} (\lambda )=0: \left\langle B_{2} , G_{m}^{0} (\lambda )\right\rangle >(1-\varepsilon)\left\langle B_{2} , U_{m} (\lambda )\right\rangle ,\]
\[\Gamma _{m} (\lambda )=0: G_{m}^{0} (\lambda )>(1-\varepsilon)U_{m}(\lambda ),\]
and $\alpha_{m} ^{2} , \beta_{m} ^{2}, \alpha _{mk}^{2} , \beta _{mk}^{2},  \vec{\alpha_{m} }, \vec{\beta_{m} }$ are the unknown Lagrange multipliers.

\begin{thm}
Let the minimality condition (\ref{eq1}) hold true. The least favourable spectral densities $F_{m}^{0}(\lambda)\in D_{0}$, $G_{m}^{0}(\lambda)\in D_{\varepsilon}$  for the optimal linear estimation of the functional  $A_{N}\zeta$ are determined by relations: (\ref{eq8}), (\ref{eq13}) -- (\ref{eq14})  for the first pair $D_{0}^{1}\times D_{\varepsilon}^{1}$,
(\ref{eq8}), (\ref{eq15}) -- (\ref{eq16})  for the second pair $D_{0}^{2}\times D_{\varepsilon}^{2}$,
(\ref{eq8}), (\ref{eq17}) -- (\ref{eq18})  for the third pair $D_{0}^{3}\times D_{\varepsilon}^{3}$,
(\ref{eq8}), (\ref{eq19}) -- (\ref{eq20})  for the fourth pair $D_{0}^{4}\times D_{\varepsilon}^{4}$.
The minimax spectral characteristic of the optimal estimate of the functional $A_{N}\zeta$ is calculated by the formula (\ref{eq2}).
\end{thm}

\begin{cor}
Let the minimality condition (\ref{eq5}) hold true. The spectral densities $F_{m}^{0}(\lambda)$ are the least favourable spectral densities in classes $D_{0}$ for the optimal linear estimation of the functional  $A_{N}\zeta$  based on observations of the field $\zeta(j,x)$ at points $j\in Z\backslash\{0, 1,\dots, N\}$,  $x\in S_{n}$ if they satisfy the following relations respectively
\begin{equation}\label{eq21}
\sum_{l=1}^{h(m,n)}((C_{mN}^{l0}(\lambda) )^{\top} )^{*}\cdot(C_{mN}^{l0}(\lambda) )^{\top}=\alpha_{m} ^{2} (F_{m}^{0} (\lambda ))^{2},
\end{equation}
\begin{equation}\label{eq22}
\sum_{l=1}^{h(m,n)}((C_{mN}^{l0}(\lambda) )^{\top} )^{*}\cdot(C_{mN}^{l0}(\lambda) )^{\top}=F_{m}^{0} (\lambda )\left\{\alpha _{mk}^{2} \delta _{kn} \right\}_{k,n=1}^{T}F_{m}^{0} (\lambda ),
\end{equation}
\begin{equation}\label{eq23}
\sum_{l=1}^{h(m,n)}((C_{mN}^{l0}(\lambda) )^{\top} )^{*}\cdot(C_{mN}^{l0}(\lambda) )^{\top}=\alpha_{m} ^{2} F_{m}^{0} (\lambda )B_{1}^{\top}F_{m}^{0} (\lambda ),
\end{equation}
\begin{equation}\label{eq24}
\sum_{l=1}^{h(m,n)}((C_{mN}^{l0}(\lambda) )^{\top} )^{*}\cdot(C_{mN}^{l0}(\lambda) )^{\top}=F_{m}^{0} (\lambda ) \vec{\alpha_{m}}\cdot \vec{\alpha_{m}}^{*}F_{m}^{0} (\lambda ),
\end{equation}
and determine solution to the extremum problem  (\ref{eq9}). The minimax spectral characteristic of the optimal estimate of the functional $A_{N}\zeta$ is calculated by the formula  (\ref{eq7}).
\end{cor}

\begin{cor}
Let the minimality condition (\ref{eq5}) hold true. The spectral densities $F_{m}^{0}(\lambda)$ are the least favourable spectral densities in classes $D_{\varepsilon}$ for the optimal linear estimation of the functional  $A_{N}\zeta$  based on observations of the field $\zeta(j,x)$ at points $j\in Z\backslash\{0, 1, \dots, N\}$,  $x\in S_{n}$ if they satisfy the following relations respectively

\begin{equation}\label{eq25}
\sum_{l=1}^{h(m,n)}((C_{mN}^{l0}(\lambda) )^{\top} )^{*}\cdot(C_{mN}^{l0}(\lambda) )^{\top}=(\alpha_{m} ^{2}+\gamma _{m} (\lambda ))(F_{m}^{0} (\lambda ))^{2},
\end{equation}

\begin{equation}\label{eq26}
\sum_{l=1}^{h(m,n)}((C_{mN}^{l0}(\lambda) )^{\top} )^{*}\cdot(C_{mN}^{l0}(\lambda) )^{\top}=F_{m}^{0} (\lambda )\left\{(\alpha _{mk}^{2}+\gamma _{mk} (\lambda ))\delta _{kn} \right\}_{k,n=1}^{T}F_{m}^{0} (\lambda ),
\end{equation}

\begin{equation}\label{eq27}
\sum_{l=1}^{h(m,n)}((C_{mN}^{l0}(\lambda) )^{\top} )^{*}\cdot(C_{mN}^{l0}(\lambda) )^{\top}=(\alpha_{m} ^{2}+\gamma_{m}^{'}(\lambda)) F_{m}^{0} (\lambda )B_{1}^{\top}F_{m}^{0} (\lambda ),
\end{equation}

\begin{equation}\label{eq28}
\sum_{l=1}^{h(m,n)}((C_{mN}^{l0}(\lambda) )^{\top} )^{*}\cdot(C_{mN}^{l0}(\lambda) )^{\top}=F_{m}^{0} (\lambda ) (\vec{\alpha_{m}}\cdot \vec{\alpha_{m}}^{*}+\Gamma _{m} (\lambda ))F_{m}^{0} (\lambda ),
\end{equation}
and determine solution to the extremum problem (\ref{eq9}). The minimax spectral characteristic of the optimal estimate of the functional $A_{N}\zeta$ is calculated by the formula (\ref{eq7}).
\end{cor}

\section*{Conclusions}

In this paper we investigate the interpolation problem for the functional
\[A_{N}\zeta=\sum _{j=0}^{N}\int_{S_{n}}a(j,x)\zeta(j,x)m_{n}(dx)\]
depending on unknown values of a spatial temporal isotropic random field
$\zeta(j,x)$, which is periodically correlated
 with respect to discrete time argument $j\in\mathrm Z$ and mean-square continuous isotropic on the sphere ${S_n}$ with respect to spatial argument $x\in{S_n}$.
 Estimates are based on
observations of the field $\zeta(j,x)+\theta(j,x)$ at
points $(j,x):$ $j\in Z\backslash\{0, 1,\dots, N\}$, $x\in{S_n},$ where $\theta(j,x)$ is an
uncorrelated with $\zeta(t,x)$
spatial temporal isotropic random field, which is periodically correlated
with respect to discrete time argument $j\in\mathrm Z$ and mean-square continuous isotropic on the sphere ${S_n}$ with respect to spatial argument $x\in{S_n}$.
The problem is investigated in the case of spectral certainty where the matrices of spectral densities of random fields are known exactly and in the case of spectral uncertainty where matrices of spectral densities of random fields are not known exactly, but some classes of admissible spectral density matrices are given.
 We derive formulas for calculation the spectral characteristic and the mean-square error of the optimal linear estimate of the functional
 provided that the matrices of spectral densities $F_{m}(\lambda ), G_{m}(\lambda )$ of the vector-valued stationary sequences that generate the random fields $\zeta(j,x)$, $\theta (j,x)$ are known exactly.
We propose a representation of the mean square
error in the form of a linear functional in the $L_1\times L_1$ space with
respect to spectral densities $(F,G)$, which allows us to solve the
corresponding conditional extremum problem and describe the minimax
(robust) estimates of the functional.
The least favourable spectral
densities and the minimax (robust) spectral characteristics of the
optimal estimates of the functional $A_N\zeta$ are determined for some
special classes of spectral densities.

\renewcommand{\refname}{References}

\end{document}